\documentclass[11pt]{article}
\usepackage{amsmath,amssymb,epsfig,setspace,framed,color}

\begin{document}

\begin{flushleft}
\textbf{\LARGE Minimal inference from incomplete $2\times 2$-tables}\\ \bigskip
Li-Chun Zhang\footnote{Primary affiliation: S3RI/Department of Social Statistics and Demography, University of Southampton, SO17 1BJ Southampton (Email: L.Zhang@soton.ac.uk). Secondary affiliation: Statistics Norway, PB. 8131 Dep, 0033 Oslo. This work was partially supported by the Administrative Data Research Centre England (ADRC-E) Grant from the UK Economic and Social Research Council (ESRC).} \\
\textit{University of Southampton, UK} \\
and Raymond L. Chambers \\
\textit{University of Wollongong, Australia}
\end{flushleft}

\noindent \textbf{Summary.} Estimates based on $2\times 2$ tables of frequencies are widely used in statistical applications. However, in many cases these tables are incomplete in the sense that the data required to compute the frequencies for a subset of the cells defining the table are unavailable. Minimal inference addresses those situations where this incompleteness leads to target parameters for these tables that are interval, rather than point, identifiable. In particular, we develop the concept of corroboration as a measure of the statistical evidence in the observed data that is not based on  likelihoods. The corroboration function identifies the parameter values that are the hardest to refute, i.e., those values which, under repeated sampling,  remain interval identified. This enables us to develop a general approach to inference from incomplete $2\times 2$ tables when the additional assumptions required to support a likelihood-based approach cannot be sustained based on the data available. This minimal inference approach then provides a foundation for further analysis that aims at making sharper inference supported by plausible external beliefs. 


\section{Introduction}

Incomplete $2\times 2$ tables are often encountered in statistical analysis. Table \ref{tab-setup} illustrates the two cases that we pay special attention to in this paper. Both tables correspond to the cross-classification of two binary variables. To the left, $X=1,0$ is the outcome variable of interest, and $R=1,0$ indicates whether an observation is missing or not. The two-way table is incomplete since $X$ is  only observed if $R=1$. We  refer to it as the \emph{missing data setting}. The two-way table on the right shows the joint distribution of two binary variables $X$ and $Y$. This table is completely unobserved. Instead, one has observations on two independent samples of $n_1$ values of $X$ and $n_2$ values of $Y$, respectively. We refer to it as the \emph{matched data setting}. For either case, we assume that the \emph{complete data} corresponding to the unobserved $2\times 2$ table follow a multinomial distribution $P_{\lambda}$, with parameter $\lambda = (\lambda_{11}, \lambda_{10}, \lambda_{01}, \lambda_{00})$ referring to the probabilities of observing each of the four possible configurations of two binary variables. Under this assumption the parameter $\lambda$ is \emph{point-identifiable} given the complete data, in the sense that $\lambda = \lambda'$ whenever $P_{\lambda} = P_{\lambda'}$; and so is the parameter of interest $\theta = \mathrm{Pr}(X=1)$ in the missing data setting and $\theta =\mathrm{Pr}(X=1, Y=1)$ in the matched data setting. 

\begin{table}[ht]
\begin{center}
\caption{Two cases of incomplete $2\times 2$-table. Left: binary variable subjected to missing data, sample size $n$; Right: statistical matching of two binary variables from separate samples of sizes $n_1$ and $n_2$, respectively. Unobserved hypothetical complete sample counts marked by `--'.}
\begin{tabular}{lcc|c||lcc|c}\hline \hline
\multicolumn{8}{c}{Hypothetical Complete Sample Data: $(n_{11}, n_{01}, n_{10}, n_{00}) \sim \text{multinomial}(n, \lambda_{11}, \lambda_{01}, \lambda_{10}, \lambda_{00})$} \\ \hline \hline
& $R=1$ & $R=0$ & Total & & $Y=1$ & $Y=0$ & Total\\ \hline
$X=1$ & $n_{11}$ & -- & -- & $X=1$ & -- & -- & $n_x$\\ 
$X=0$ & $n_{01}$ & -- & -- & $X=0$ & -- & -- & $n_1 - n_x$\\ \hline
Total & $n_{+1}$ & $n_{+0}$ & $n$ & Total & $n_y$ & $n_2 - n_y$ &  \\ \hline 
Observation: & \multicolumn{3}{l||}{$(n_{11}, n_{01}, n_{+0})$ Given $n$} & 
Observation: & \multicolumn{3}{l}{ Independent $(n_x, n_1)$, $(n_y, n_2)$} \\ \hline
Sampling & \multicolumn{3}{c||}{} &
Sampling & \multicolumn{3}{l}{$n_x\sim \text{binomial}(n_1, \lambda_{1+})$}  \\
Distribution: & \multicolumn{3}{l||}{multinomial$(n, \lambda_{11}, \lambda_{01}, \lambda_{+0})$} & 
Distribution: & \multicolumn{3}{l}{$n_y\sim \text{binomial}(n_2, \lambda_{+1})$} \\ \hline 
Identifiable: & \multicolumn{3}{l||}{$\lambda_{11}, \lambda_{01}, \lambda_{+0} = \lambda_{10} + \lambda_{00}$} &
\multicolumn{4}{l}{Identifiable: $\lambda_{1+} = \lambda_{11} +\lambda_{10}, \lambda_{+1} = \lambda_{11} + \lambda_{01}$} \\ \hline
\multicolumn{4}{l||}{Parameter of Interest: $\theta = \lambda_{1+}$} & 
\multicolumn{4}{l}{Parameter of Interest: $\theta = \lambda_{11}$}\\ \hline
\multicolumn{4}{l||}{Identification Region: $\lambda_{11}\leq \theta \leq \lambda_{11}+\lambda_{+0}$} &
\multicolumn{4}{l}{Identification Region: $\theta \leq \min(\lambda_{1+}, \lambda_{+1})$}\\
& \multicolumn{3}{c||}{} & & \multicolumn{3}{r}{$\theta \geq \max(\lambda_{1+} + \lambda_{+1}-1, 0)$} \\ \hline
\multicolumn{4}{l||}{Additional Assumption: Independent $(X,R)$} & 
\multicolumn{4}{l}{Additional Assumption: Independent $(X,Y)$} \\ \hline \hline
\end{tabular}\label{tab-setup} \end{center} \end{table} 

Table \ref{tab-setup} also shows the sampling distribution of the observed data for each setting. Point-identification for $\theta$ based on the observed data is only achievable if additional assumptions are made. In this table these assumptions are independence of $X$ and $R$ in the missing data setting, which means missing-completely-at-random (MCAR, Rubin, 1976), and independence of $X$ and $Y$ in the statistical matching setting, which is a special case of the conditional independence assumption (Okner, 1972). But such additional assumptions are often contentious. It therefore seems  reasonable ask `what the data say' about $\theta$ given the \emph{accepted} sampling distribution of the observed data, \emph{without} the additional ``esoteric'' (Tamer, 2010) assumptions that enable point-identification of this parameter. The aim of this paper is to describe a general approach to inference based on incomplete $2\times 2$ tables given such a setting. 

\begin{figure}[ht]
\resizebox{16cm}{8cm}{\includegraphics{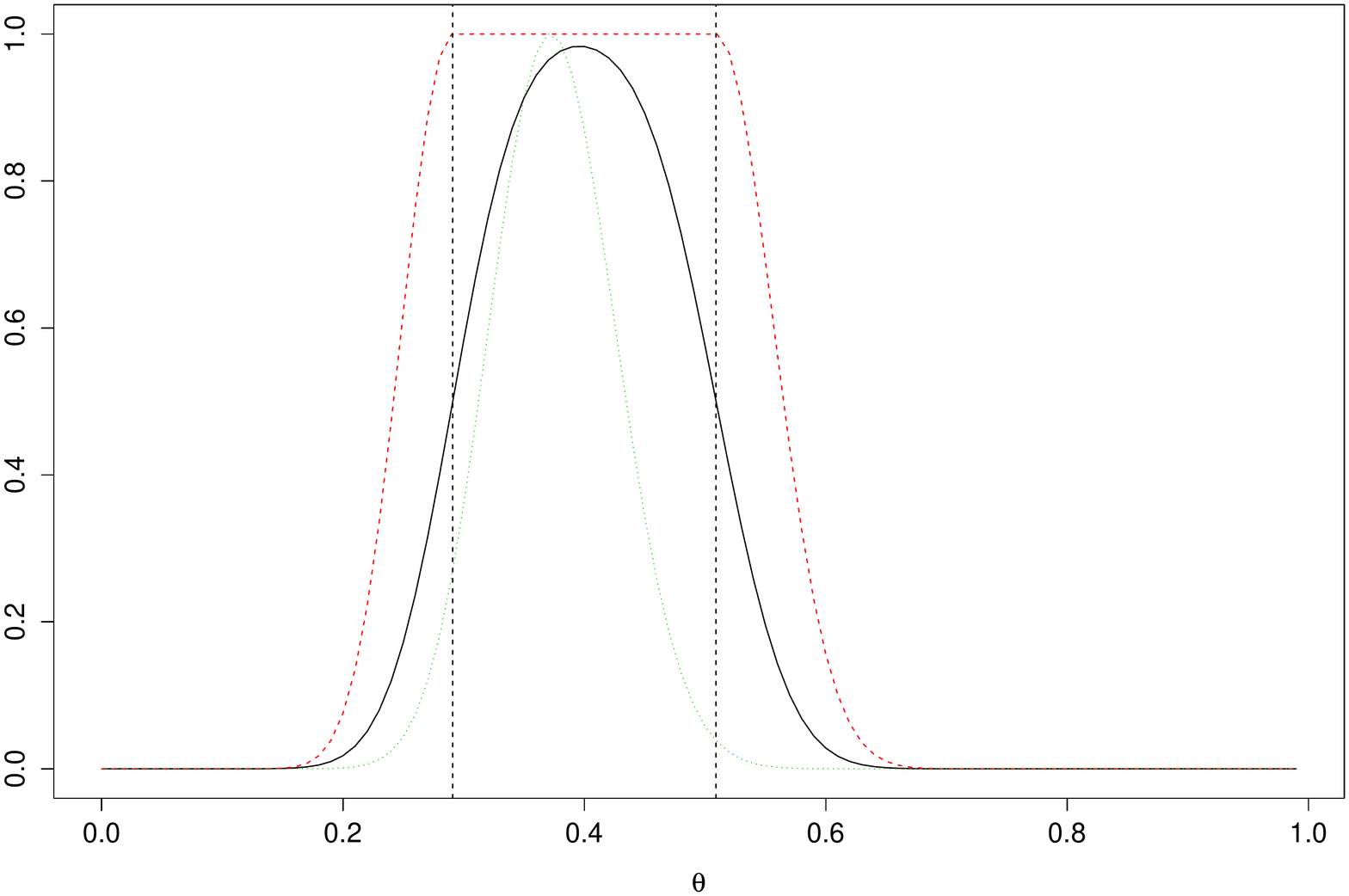}}
\caption{Observed corroboration and standardised likelihoods (with peak value $1$) based on OCBGT data: observed corroboration (solid), profile likelihood (dashed), likelihood under MCAR assumption (dotted); maximum likelihood region marked by vertical dashed lines.} \label{fig-OCBGT} 
\end{figure}

To illustrate, consider a missing data example discussed by Zhang (2010). The observed data from the Obstructed Coronary Bypass Graft Trials (OCBGT, see Hollis, 2002) are $(n_{11}, n_{01}, n_{+0}) = (32, 54, 24)$, with the sampling distribution parameter $\psi = (\lambda_{11},\lambda_{01},\lambda_{+0})$. The likelihood of $\psi$ is proportional to $\lambda_{11}^{n_{11}} \lambda_{01}^{n_{01}} \lambda_{+0}^{n_{+0}}$. This yields the profile likelihood of the parameter of interest $\theta = \lambda_{1+}$, denoted by $L_p(\theta)$, which is the dashed curve in Figure \ref{fig-OCBGT}. It is seen that $L_p(\theta)$ is flat over $[n_{11}/n, (n_{11} + n_{+0})/n]$, which we call the \emph{maximum likelihood region}, denoted by $\widehat{\Theta}$, with all values of $\theta$ in $\widehat{\Theta}$ equally likely based on the observed data. Asymptotically, as $n\rightarrow \infty$, $\widehat{\Theta}$ tends to the \emph{identification region} of $\theta$, i.e. $\lambda_{11}\leq \theta \leq \lambda_{11}+\lambda_{+0}$, which is a function of the identifiable parameter $\psi$. This identification region is the asymptote of `what the data say' about $\theta$ under the setting here. The dotted curve gives the standardised likelihood under the additional MCAR assumption that enables point-identification of $\theta$. It peaks at the maximum likelihood estimate (MLE) $\widehat{\theta}_{MCAR} = n_{11}/n_{+1}$, which converges to $\lambda_{11}/\lambda_{+1}$ in probability. Clearly, the MLE derived from the MCAR likelihood will be \emph{inconsistent} as long as $\lambda_{1+} \ne \lambda_{11}/\lambda_{+1}$.

The fact that the profile likelihood shown in Figure \ref{fig-OCBGT} is constant within the observed $\widehat{\Theta}$ does not mean that all the values of $\theta$ in it are equally likely to be in a $\widehat{\Theta}$ that could be observed given a random draw from the sampling distribution of the observed data. In Section \ref{definition} we develop the concept of \emph{corroboration}, noting that values of $\theta$ that are more likely to appear in a $\widehat{\Theta}$ on repeated sampling are better corroborated by the observed data than values of $\theta$ that only infrequently appear in a $\widehat{\Theta}$. The solid curve in Figure \ref{fig-OCBGT} shows how the estimated corroboration varies with $\theta$ for the OCBGT data. The computation of the estimated corroboration is explained in Section \ref{definition}. The key point to note here is that the corroboration varies for the points within $\widehat{\Theta}$, where the profile likelihood is constant. This allows us to construct high corroboration level sets \emph{within} $\widehat{\Theta}$. It will be shown that asymptotically the set of values with the maximum observed corroboration becomes indistinguishable from the identification region except for its bounds. Unlike the MLE that aims at the \emph{most likely} parameter value, the maximum corroboration set identifies those parameter values that are the \emph{hardest to refute} based on the observed data. In effect, these are the points in which we have the highest confidence. We develop a Corroboration Test in Section \ref{ctest} for the settings of Table \ref{tab-setup}, where the Likelihood Ratio Test is inapplicable insofar as the parameter of interest is not point identifiable. The test will be applied the OCBGT data.

There are several related approaches within the matched data setting. In ecological inference (Goodman, 1953; King, 1997), the observed data are the margins of the unobserved complete $2 \times 2$ table. See Wakefield (2004) for a comprehensive review. It is clearly recognised that critical but untestable assumptions are needed to arrive at a point estimate in this context, and that there is a fundamental difficulty associated with choosing between different models based only on the observed data; see Greenland and Robins (1994), Freedman (2001) and Gelman \emph{et al.} (2001). Statistical matching deals with the same setting, where the set of multinomial distributions $P_{\lambda}$ compatible with the sampling distribution of the observed data is referred to as the \emph{uncertainty space}. Evaluation of the uncertainty space has received much attention (Kadane, 1978; Moriarity and Scheuren, 2001; D'Orazio {\em et al.}, 2006; Kiesel and R\"{a}ssler, 2006; Conti {\em et al.} 2012; Zhang, 2015; Conti \emph{et al.}, 2015). The concept of uncertainty space is closely related to that of identification uncertainty (Koopmans, 1949; Tamer, 2010). The ``partial identification'' framework (Manski, 1995, 2003, 2007) recognises situations where, due to the structure of the data, even a hypothetical infinite number of observations may only constrain the parameter of interest without being able to point-identify it. It is important in this context to distinguish between the study of identification, provided an \emph{infinite} amount of data under the given structure, and statistical inference from \emph{finite} samples. Partial identification in econometrics can be traced back to Frisch (1934) and Marschak and Andrews (1944), and there is a growing literature on the construction of confidence regions of the identified parameter set. See e.g. Imbens and Manski (2004), Chernozhukov \emph{et al.} (2007), Beresteanu and Molinari (2008), and Ramano and Shaikh (2010).

It is clear that all the aforementioned approaches aim at inference based on an identifiable sampling distribution that is acceptable to all, no matter which untestable additional assumptions an analyst may or may not introduce in order to resolve the identification issue. As seen in Figure \ref{fig-OCBGT}, the novelty of the approach proposed in this paper is that it achieves this objective via a measure of the statistical evidence in the observed data that is \emph{not} based on comparing likelihoods.

\section{Corroboration} \label{definition}

Denote by $f(d_n;\psi)$ the identifiable sampling distribution of the observed data $d_n$ with generic sample size $n$, and with parameter $\psi$. Denote by $P_{\lambda}$ the distribution of the hypothetical complete data, which is characterised by the parameter $\lambda$ with parameter space $\Lambda$. Denote by $\theta = \theta(\lambda)$ a scalar parameter of interest, and by $\Theta$ the parameter space of $\theta$. For any given $\psi$ let $\Lambda(\psi)$ be the constrained parameter space defined by $\psi$. That is, $\Lambda(\psi)$ consists of all $\lambda$ that are consistent with $\psi$. Let $\Theta(\psi)$ be the induced parameter space of $\theta$, which contains all $\theta(\lambda)$ where $\lambda \in \Lambda(\psi)$. For inference under a \emph{minimal setting} in this paper, we then require both conditions below to hold.
\begin{description}
\item[(M$_1$)] The induced parameter space $\Theta(\psi)$ is a closed interval. In particular, it is not a singleton $\Theta(\psi) = \theta(\psi)$, nor is it invariant towards $\psi$ in the sense that $\Theta(\psi) = \Theta(\psi')$ for all $\psi \neq \psi'$.

\item[(M$_2$)] The parameter $\psi$ of the sampling distribution is point-identifiable, and the MLE $\widehat{\psi}$ is such that $\widehat{\psi} \stackrel{\mathrm{Pr}}{\rightarrow} \psi_0$, asymptotically as $n\rightarrow \infty$, where $\psi_0$ is the true parameter value.
\end{description}
Under a minimal setting, $\Theta(\psi) = [L(\psi), U(\psi)]$, where $L(\psi)$ is the \emph{lower bound} of $\theta$ induced by $\psi$, and $U(\psi)$ the \emph{upper bound}. The identification region is $\Theta_0 = \Theta(\psi_0) = [L_0, U_0]$, where $L_0 = L(\psi_0)$ and $U_0 = U(\psi_0)$. Thus, for the missing data setting in Table \ref{tab-setup}, we have $\psi_0 = (\lambda_{11}^0, \lambda_{01}^0, \lambda_{+0}^0)$, with
\[
\Theta_0 = [L_0, U_0] = [\lambda_{11}^0,~  \lambda_{11}^0 + \lambda_{+0}^0]. 
\] 
For the matched data setting, we have $\psi_0 = (\lambda_{1+}^0, \lambda_{+1}^0)$, and the Fr\'{e}chet bounds (Fr\'{e}chet, 1951) define the identification region
\[
\Theta_0 = [L_0, U_0] = [ \max(\lambda_{1+}^0 + \lambda_{+1}^0 -1, 0),~  \min(\lambda_{1+}^0, \lambda_{+1}^0)]. 
\]

Let $\widehat{L} = L(\widehat{\psi})$ and $\widehat{U} = U(\widehat{\psi})$ be the MLEs of $L_0$ and $U_0$, respectively, and let  $\widehat{\Theta} = \Theta(\widehat{\psi}) = [\widehat{L}, \widehat{U}]$ denote the maximum profile likelihood estimator of $\theta$. The points inside $\widehat{\Theta}$ can all be considered as \emph{equally most likely}, i.e. best supported according to the likelihood based on $d_n$ under the observed data model. We define the \emph{corroboration function} of $\theta$, for $\theta \in \Theta$, to be
\begin{equation}
c(\theta; \psi) = \mathrm{Pr}(\theta \in \widehat{\Theta}; \psi),  \label{c-func}
\end{equation}
i.e. the probability for the given value of $\theta$ to be covered by $\widehat{\Theta}$, where the probability is evaluated with respect to $f(d_n; \psi)$. Let the \emph{actual corroboration} be 
\[
c_{0}(\theta) = c(\theta; \psi_0),
\]
i.e. evaluated over the true sampling distribution. In particular, $c(\theta_0; \psi_0)$ is the confidence level of $\widehat{\Theta}$ as an interval estimator of $\theta_0$. Let the \emph{observed corroboration} be
\[
\widehat{c}(\theta) = c(\theta; \widehat{\psi}). 
\]
Since $\widehat{c}(\theta)$ is the MLE of $c_{0}(\theta)$, one may then define the observed corroboration as the \textit{most likely level of corroboration} for $\theta$ given the observed data. As illustrated in Figure \ref{fig-OCBGT} for the OCBGT data, if one treats the observed corroboration as a function of $\theta$ then this function can generally vary over $\widehat{\Theta}$, as opposed to the profile likelihood which is flat over the same region. Note that in this case in order to calculate $\widehat{c}(\theta)$, where $(\widehat{\lambda}_{11}, \widehat{\lambda}_{+0}) = (n_{11}/n, n_{11}/n + n_{+0}/n)$, we employ the bivariate normal approximation $(\widehat{\lambda}_{11}, \widehat{\lambda}_{+0}) \sim N_2(\mu, \Sigma)$, where $\mu = (\lambda_{11}, \lambda_{+0})$ and the distinctive elements of $\Sigma$ are $V(\widehat{\lambda}_{11}) = \lambda_{11} (1- \lambda_{11})/n$, $V(\widehat{\lambda}_{+0}) = \lambda_{+0} (1- \lambda_{+0})/n$ and $Cov(\widehat{\lambda}_{11}, \widehat{\lambda}_{+0}) = - \lambda_{11} \lambda_{+0}/n$. More generally, the observed corroboration can be calculated via simulation as follows. 

\paragraph{Bootstrap for $\widehat{c}(\theta)$} For given $\theta$ and the MLE $\widehat{\psi}$, repeat for $b=1, ... B$:
\begin{itemize}
\item generate $d_n^{(b)}$ from $f(d_n; \widehat{\psi})$ to obtain $\widehat{\psi}^{(b)}$ and the corresponding $[L(\widehat{\psi}^{(b)}), U(\widehat{\psi}^{(b)})]$; 

\item set $\delta^{(b)} = 1$ if $\theta \in [L(\widehat{\psi}^{(b)}), U(\widehat{\psi}^{(b)})]$, and 0 otherwise.    
\end{itemize} 
Put $\widehat{c}({\theta}) = \sum_{b=1}^B \delta^{(b)}/B$ as the bootstrap estimate of the observed corroboration for $\theta$. $\square$

\section{Maximum corroboration set} \label{MACS}

Let the \emph{level-$\alpha$ corroboration set} be given by  
\[
A_{\alpha}(\psi) = \{ \theta : c(\theta; \psi) \geq \alpha \}, 
\]
\emph{provided} there exists some $\theta \in A_{\alpha}(\psi)$ where $c(\theta; \psi) = \alpha$. Thus, by definition we have $c(\theta; \psi) < \alpha$, for any $\theta\not \in A_{\alpha}(\psi)$, whilst we \emph{cannot} have $c(\theta; \psi) > \alpha$ for all $\theta\in A_{\alpha}(\psi)$. Some properties of $A_{\alpha}(\psi)$ are given below, with proofs in the Appendix. Notice that we use $c(\theta)$ as a short-hand for $c(\theta; \psi)$ and $A_{\alpha}$ that of $A_{\alpha}(\psi)$, where it is not necessary to emphasise their dependence on $\psi$.

\paragraph{Theorem 1} Suppose that a minimal inference setting applies, i.e. provided conditions (M$_1$) and (M$_2$) hold. Then:
\begin{description}
\item[(i)] Let $A_{\alpha_1} = [L_1, U_1]$ and $A_{\alpha_2} = [L_2, U_2]$. If $\alpha_1 > \alpha_2$, then $[L_1, U_1] \subset [L_2, U_2]$.

\item[(ii)]  Let $\theta_L < \theta_U$, where $c(\theta_L) = c(\theta_U) =\alpha$. Then $c(\theta) \geq \alpha$ for any $\theta \in (\theta_L, \theta_U)$. 
\end{description}

\paragraph{Theorem 2} Given a minimal inference setting, there exists a \emph{maximum corroboration value} denoted by $\theta^{\max}$, such that $c(\theta^{\max}) \geq c(\theta)$ for any $\theta \neq \theta^{\max}$.

\bigskip
Denote by $A^{\max} = A^{\max}(\psi_0)$ the \emph{maximum corroboration set}, such that $c_{0}(\theta) > c_{0}(\theta')$ for any $\theta \in A^{\max}$ and $\theta' \not \in A^{\max}$, and $c_{0}(\theta) = c_{0}(\theta')$ for any $\theta \neq \theta' \in A^{\max}$. It follows from \eqref{c-func} that these are the points for which $\widehat{\Theta}$ implies the highest confidence, in which sense one may consider these to be the parameter values that are the hardest to refute. Replacing $\psi_0$ by $\widehat{\psi}$, we obtain the MLE of $A^{\max}$ or the \emph{observed} maximum corroboration set
\[
\widehat{A}^{\max} = A^{\max}(\widehat{\psi}).
\]

\begin{figure}[ht]
\resizebox{16cm}{7cm}{\includegraphics{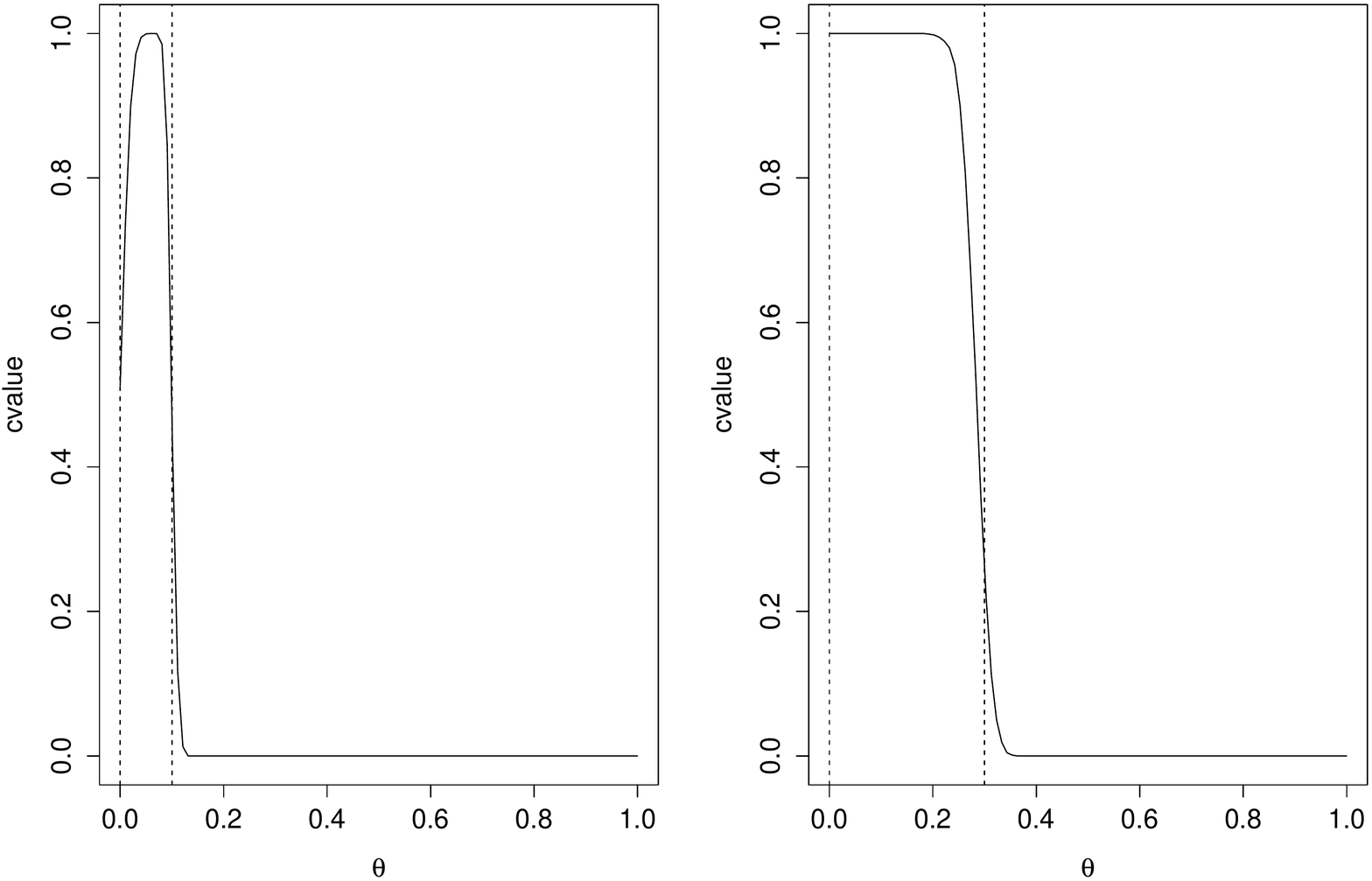}}
\caption{Illustration of corroboration in matched data setting. Left: $(\lambda_{1+}, n_1) = (0.1, 1000)$ and $(\lambda_{+1}, n_2) = (0.9, 500)$. Right: $(\lambda_{1+}, n_1) = (0.3, 200)$ and $(\lambda_{+1}, n_2) = (0.3, 300)$.} \label{fig-asymp} 
\end{figure}

Figure \ref{fig-asymp} illustrates corroboration in the matched data setting, where $\theta = \lambda_{11}$. The true sampling distribution parameters $(\lambda_{1+}, \lambda_{+1})$ are $(0.1,0.9)$ for the left plot and $(0.3, 0.3)$ to the right. The sample sizes are $(n_1, n_2) = (1000, 500)$ to the left and $(200, 300)$ to the right. The identification region $\Theta_0$ is the interval between the vertical dashed lines, and the solid curve shows how the actual corroboration (denoted $cvalue$ in the plots) varies with $\theta$. The corroboration of some interior points of $\Theta_0$ can be 1, whereas it can be 0 for many $\theta \not \in \Theta_0$. In the left plot, both $c_{0}(L_0)$ and $c_{0}(U_0)$ are about 0.5; in the right plot, we have $c_{0}(L_0) = 1$ and $c_{0}(U_0) \approx 0.25$. 

\begin{table}[ht]
\begin{center}
\caption{Asymptotic actual corroboration $\bar{c}_{0}(\theta)$ in missing and matched data settings}
\begin{tabular}{|l|c|c|c|c|}\hline \hline
Data Setting & $\theta \not \in [L_0, U_0]$ & $\theta = L_0$ & $\theta \in (L_0, U_0)$ & $\theta = U_0$ \\ \hline
Missing  & 0 & 0.5 if $L_0 >0$ & 1 & 0.5 if $U_0 <1$\\ \hline
Matching  & 0 & 0.5 if $\lambda_{1+} + \lambda_{+1} \geq 1$ & 1 & 0.5 if $\lambda_{1+} \neq \lambda_{+1}$ \\  
 & & 1 if $\lambda_{1+} + \lambda_{+1} < 1$ & & 0.25 if $\lambda_{1+} = \lambda_{+1}$ \\ \hline \hline
\end{tabular}\label{tab-asymptote} \end{center} \end{table}

Let $\bar{c}(\theta; \psi) = \lim_n c(\theta; \psi) = \lim_n \mathrm{Pr}(\theta \in \widehat{\Theta}_n; \psi)$ be the \emph{asymptotic corroboration} of $\theta$ evaluated at $\psi$,  where $\lim_n$ stands for $\lim_{n\rightarrow \infty}$ and $\widehat{\Theta}_n$ makes explicit the dependence on sample size. Table \ref{tab-asymptote} summarises the asymptotic actual corroboration $\bar{c}_{0}(\theta) = \bar{c}(\theta; \psi_0)$ for both data settings. Let $\bar{A}^{\max}$ be the asymptotic maximum actual corroboration set based on $\bar{c}_{0}(\theta)$. Lemma 1 states that, apart from the bounds $L_0$ and $U_0$, $\bar{A}^{\max}$ is indistinguishable from $\Theta_0$ and $\bar{c}_0(\theta)$ is an indicator function on $\Theta_0$. Theorem 3 states that the interior of the observed maximum corroboration set $\widehat{A}_n^{\max}$ converges to the interior of $\Theta_0$ in probability. The proofs are given in the Appendix.

\paragraph{Lemma 1} Given a minimal inference setting, $\theta \in \bar{A}^{\max}$ and $\bar{c}_0(\theta) =1$ if $\theta \in \mbox{Int}(\Theta_0) = (L_0, U_0)$, i.e. if $\theta$ belongs to the interior of $\Theta_0$, then $\theta \not \in \bar{A}^{\max}$ and $\bar{c}_0(\theta) =0$ for any $\theta \not \in [L_0, U_0]$.

\paragraph{Theorem 3} Given a minimal inference setting, we have $\mbox{Int}(\widehat{A}^{\max}) \stackrel{\mbox{Pr}}{\rightarrow} \mbox{Int}(\Theta_0)$; that is, $\lim_n \mbox{Pr}(\theta \in \widehat{A}_n^{\max}) =1$ if $\theta \in \mbox{Int}(\Theta_0)$ and $\lim_n \mbox{Pr}(\theta \in \widehat{A}_n^{\max}) =0$ if $\theta \not \in \Theta_0$.

\section{High assurance estimation of $\Theta_0$} \label{assurance}

Given a minimal inference setting, a confidence region $C_n$ for $\Theta_0$ (which is an interval) has the confidence level $\mbox{Pr}(\Theta_0 \subseteq C_n)$; see e.g. Chernozhukov \emph{et al.} (2007). Given a high confidence level, the probability that $C_n$ contains points that do not belong to $\Theta_0$ must also be high, due to sampling variability, and so $C_n$ asymptotically contracts towards $\Theta_0$ from `outside' of it. In contrast, any point in $\Theta_0$ is irrefutable, and $\widehat{A}^{\max}$ identifies those parameter values that are the hardest to refute given the observed data. We thus define the \emph{assurance} of $\widehat{A}^{\max}$ to be 
\[
\tau_0 = \mathrm{Pr}(\widehat{A}^{\max} \subseteq \Theta_0), 
\]
where the probability is evaluated with respect to $f(d_n; \psi_0)$. That is, this is the probability that the points in the observed $\widehat{A}^{\max}$ are indeed all irrefutable. If $\widehat{A}^{\max}$ has a high assurance, there will be a low probability that it contains points outside of $\Theta_0$. As the sample size increases, a high assurance estimator of $\Theta_0$ should therefore grow towards $\Theta_0$ from `inside' of it. In light of Theorem 1, for some small constant $h\geq 0$, a high assurance estimator of $\Theta_0$ can therefore be defined as
\[
\widehat{A}_h = \{ \theta: c(\theta; \widehat{\psi}) \geq \max_{\theta} c(\theta; \widehat{\psi}) - h \},
\] 
The following bootstrap can be used to estimate $\widehat{A}_h$, including $\widehat{A}_0 = \widehat{A}^{\max}$.

\paragraph{Bootstrap for $\widehat{A}_h$} Given the MLE $\widehat{\psi}$ and the corresponding $[\widehat{L}, \widehat{U}]$, repeat for $b=1, ... B$:
\begin{enumerate}
\item generate $d_n^{(b)}$ from $f(d_n; \widehat{\psi})$, and obtain $\widehat{\psi}^{(b)}$; 

\item for any given $h$, where $0\leq h<1$, obtain $\widehat{A}_h^{(b)}$ at $\widehat{\psi}^{(b)}$ in the same way as $\widehat{A}_h$ at $\widehat{\psi}$, and the corresponding $L^{(b)} = L(\widehat{A}_h^{(b)})$ and $U^{(b)} = U(\widehat{A}_h^{(b)})$;

\item set $\delta^{(b)} = 1$ if $\widehat{L} \leq L^{(b)}  <  U^{(b)} \leq \widehat{U}$, and $\delta^{(b)} = 0$ otherwise.    
\end{enumerate} 
Calculate the bootstrap estimate of assurance as $\widehat{\tau}(\widehat{A}_h; \psi_0) = \sum_{b=1}^B \delta^{(b)}/B$, with corresponding bootstrap estimate of the lower end of $\Theta_0$ given by $L(\widehat{A}_h) = \sum_{b=1}^B L^{(b)}/B$ and of the upper end of $\Theta_0$ given by $U(\widehat{A}_h) = \sum_{b=1}^B U^{(b)}/B$. $\square$ 

\bigskip
For small $h$, $\widehat{A}_h$ can have higher assurance than $\widehat{\Theta}$, whereas it can be `closer' to $\bar{A}^{\max}$ than $\widehat{A}^{\max} = \widehat{A}_0$ by Theorem 1, since $\widehat{A}_0 \subset \widehat{A}_h$. Setting $h< 0.25$ makes $\mbox{Int}(\widehat{A}_h)$ asymptotically indistinguishable from $\mbox{Int}(\Theta_0)$ for the two settings depicted in Table \ref{tab-asymptote}. In a finite-sample situation, one may calculate $\widehat{A}_h$ and its assurance for several different choices of $h$. Since the length of $\widehat{A}_h$ increases with $h$ while its assurance decreases, one may choose the longest $\widehat{A}_h$ as an estimator of $\Theta_0$ subject to an acceptable level of assurance.

\section{A Corroboration Test} \label{ctest}

Consider testing the null hypothesis $H_A: \theta^* \in (L_0, U_0)$ against $H_B: \theta^* \not \in \Theta_0$. A minimal inference setting for this test is nonstandard because, under both $H_A$ and $H_B$, the set of possible distributions of the observed data are exactly the same, i.e. $f(d_n;\psi)$. The Likelihood Ratio Test is inapplicable. Let instead the test statistic be $T_n =1$ if $\theta^* \in \mbox{Int}(\widehat{\Theta}_n)$ and $T_n =0$ if $\theta^* \not \in \widehat{\Theta}_n$. Suppose we reject $H_A$ if $T_n =0$. The power function of this testing procedure is then $\beta_n(\theta^*) = \mbox{Pr}(T_n = 0; \psi_0)$, and is such that
\[
\bar{\beta}(\theta^*) \equiv \lim_n \beta_n(\theta^*) = 1 - \lim_n \mbox{Pr}(T_n =1; \psi_0) = 1- \bar{c}_0(\theta^*). 
\]
If $H_A$  is true, but $T_0 =0$ and we reject $H_A$, by Lemma 1 the probability of Type-I error converges to zero since $\bar{c}_0(\theta^*) = 1$ if $\theta \in \mbox{Int}(\Theta_0)$. Similarly, if $H_B$ is true, but $T=1$ and we do not reject $H_A$, the Type-II error probability also asymptotes to zero since $\bar{c}_0(\theta^*) = 0$ if $\theta^* \not \in \Theta_0$.

\begin{table}[ht]
\centering
\caption{Supporting evidence for $H_A: \theta^* \in (L_0, U_0)$ vs. $H_B: \theta^* \not \in \Theta_0$.}
\begin{tabular}{|c|c|c|} \hline \hline
& Low Power $\widehat{\beta}_n(\theta^*)$ & High Power $\widehat{\beta}_n(\theta^*)$ \\ \hline
$T_n =1$ & Support $H_A$ & Support neither, improbable event \\ \hline
$T_n =0$ & Support neither, improbable event & Support $H_B$ \\ \hline\hline
 \end{tabular}\label{tab-test}  \end{table} 

Let the \emph{observed power} be $\widehat{\beta}_n(\theta^*) = 1 -\widehat{c}_n(\theta^*)$, which is a consistent estimator of $\bar{\beta}(\theta^*)$. While $\widehat{c}_n(\theta^*)$ is a consistent estimator of the Type-II error probability, we cannot use it to estimate the Type-I error probability. The reason is that $c_0(\theta^*)$ is the same under $H_A$ or $H_B$, due to the minimal inference setting, so that it cannot be related to \emph{both} types of errors. We shall therefore define the \emph{Corroboration Test} to have observed power $\beta$, where $\beta = \widehat{\beta}_n(\theta^*) \in (0,1)$, if $H_A$ is rejected when $T_n = 0$. As summarised in Table \ref{tab-test}, a Corroboration Test of high observed power would lead one to reject $\theta^*$ if it is outside of $\widehat{\Theta}_n$ and have a low observed corroboration. By the consistency of $\widehat{A}_n^{\max}$ established in Theorem 3, we have 
\[
\lim_n \mathrm{Pr}(\text{Reject } H_A \text{ when } H_A \text{ is true}) = 0 < \lim_n \mathrm{Pr}(\text{Reject } H_A \text{ when } H_B \text{ is true}) = 1. 
\]
That is, the Corroboration Test is strongly Chernoff-consistent, since $T_n$ has limiting size 0 and the Type-II error probability converges to 0, for any $\theta^*$ specified in $H_A$.

\paragraph{Theorem 4} Given a minimal inference setting, the Corroboration Test of observed power $\beta = \widehat{\beta}_n(\theta^*)$, for $\beta \in (0,1)$, is strongly Chernoff-consistent.

\section{Application: Missing OCBGT data} \label{application}

Consider the OCBGT data $n = (n_{11}, n_{01}, n_{+0}) = (32, 54, 24)$. The profile likelihood is  
\[
L_p(\theta) \propto \begin{cases}�
n_{11}^{\theta} n_{01}^{\frac{(1-\theta) n_{01}}{n_{01} + n_{+0}}} n_{+0}^{\frac{(1-\theta) n_{+0}}{n_{01} + n_{+0}}} & \text{if } \theta < \widehat{\lambda}_{11} \\
n_{11}^{\widehat{\lambda}_{11}} n_{01}^{\widehat{\lambda}_{01}} n_{+0}^{\widehat{\lambda}_{+0}} & \text{if } \widehat{\lambda}_{11} \leq \theta \leq \widehat{\lambda}_{11} + \widehat{\lambda}_{+0} \\
n_{11}^{\frac{\theta n_{11}}{n_{11} + n_{+0}}} n_{01}^{1-\theta} n_{+0}^{\frac{\theta n_{+0}}{n_{11} + n_{+0}}} & \text{if } \theta > \widehat{\lambda}_{11} + \widehat{\lambda}_{+0}
\end{cases}
\] 
(Zhang, 2010). The likelihood is $L_{MCAR}(\theta) \propto n_{11}^{\theta} n_{01}^{1-\theta}$, under the additional assumption of independent $(X,R)$. Figure \ref{fig-OCBGT} plots both, as well as the observed corroboration $\widehat{c}(\theta)$. 

The likelihood $L_{MCAR}$ does not vary with $n_{+0}$, e.g. whether this value is 4, 24 or 104. Accordingly $n_{+0}$ is not part of the available statistical evidence. Clearly, such insensitiveness towards the observed data requires some external belief to sustain. Next, consider the relative plausibility of $\theta^* = 0.2, 0.3, 0.5, 0.6$ against $\theta_1 = 0.4$ based on the profile likelihood ratio, denoted by $LR_p(\theta^*, \theta_1)$ in the left part of Table \ref{tab-OCBGT}. The values $0.3$ and $0.5$ cannot be distinguished from $0.4$, since all are inside $\widehat{\Theta} = [0.29, 0.51]$; the negative evidence of $0.2$ and $0.6$ against $0.4$ is ``moderate'' according to Royall (1997), as they fall in the range $1/32 - 1/8$. Nevertheless, as noted before, the Likelihood Ratio Test is inapplicable here.

\begin{table}[ht]
\centering
\caption{Left, profile likelihood ratio $LR_p(\theta^*,\theta_1)$ with $\theta_1 =0.4$, observed corroboration $\widehat{c}(\theta^*)$ based on OCBGT data. Right, assurance $\widehat{\tau}(\widehat{A}_h; \psi_0)$ of $\widehat{A}_h$, expected left end $L(\widehat{A}_h)$ and right end $U(\widehat{A}_h)$, with values obtained by bootstrap with  $B=5000$. In addition, $\widehat{\Theta} : [\widehat{L}, \widehat{U}] = [0.29, 0.51]$,  $\widehat{\tau}(\widehat{\Theta}) = 0.19$.} 
\begin{tabular}{ccc||ccc} \hline \hline
$\theta^*$ & $LR_p(\theta^*, \theta_1)$ & $\widehat{c}(\theta^*)$ & 
$h$ & $\widehat{\tau}(\widehat{A}_h; \psi_0)$ & $[L(\widehat{A}_h), U(\widehat{A}_h)]$ \\ \hline
0.2 & 0.076 & 0.018 & 0 & 0.99 & [0.40, 0.40] \\
0.3 & 1 & 0.583 & 0.01 & 0.95 & [0.38, 0.41] \\
0.4 & 1 & 0.985 & 0.06 & 0.84 & [0.36, 0.44] \\
0.5 & 1 & 0.576 & 0.40 & 0.25 & [0.30, 0.50] \\
0.6 & 0.156 & 0.028 & 0.80 & 0.00 & [0.25, 0.55] \\ \hline \hline
 \end{tabular}\label{tab-OCBGT} \end{table} 

Now, based on the observed corroboration $\widehat{c}(\theta^*)$ in Table \ref{tab-OCBGT}, one may reject the null hypothesis $H_0: 0.2 \in \Theta_0$ on the basis of the Corroboration Test with observed power $0.982$. Similarly for $H_0: 0.6\in \Theta_0$, with observed power $0.972$. Meanwhile, $0.3$ and $0.5$ are just inside $\widehat{\Theta}$, with $\widehat{c}(0.3)$ and $\widehat{c}(0.5)$ slightly below 0.6, and so cannot be rejected with high observed power. The Corroboration Test thus allows us to reject an unlikely value of $\theta$ with a high observed power.  

Finally, five observed corroboration level sets $\widehat{A}_h$ are illustrated in the right part of Table \ref{tab-OCBGT}, where the estimated assurance $\widehat{\tau}(\widehat{A}_h; \psi_0)$ and expected end points $L(\widehat{A}_h)$ and $U(\widehat{A}_h)$ are calculated using the bootstrap described in Section \ref{assurance}. As an estimator of $\Theta_0$, $\widehat{A}_0$ is very narrow but has $99\%$ assurance; $\widehat{A}_{0.01}$ has $95\%$ assurance and is expected to span from $0.38$ to $0.41$. Using $\widehat{\Theta}$ as an estimator of $\Theta_0$ would perform comparably to $\widehat{A}_{0.4}$, but with low assurance. The observed corroboration level sets thus allow us to identify true irrefutable points in $\Theta_0$ with a high assurance.

\appendix

\section{Appendix}
\subsection{Proof of Theorem 1}

\noindent
(i) On the one hand, we have $A_{\alpha_1} \setminus A_{\alpha_2} = \emptyset$ because, otherwise, there must exist some $\theta \in A_{\alpha_1} \setminus A_{\alpha_2}$ such that $c(\theta) \geq \alpha_1$ (because $\theta \in A_{\alpha_1}$) and $c(\theta) < \alpha_2$ (because $\theta \not \in A_{\alpha_2}$) at the same time, contradictory to $\alpha_1 >\alpha_2$ as stipulated. On the other hand, the set $A_{\alpha_2} \setminus A_{\alpha_1}$ is non-empty because, otherwise, every $\theta\in A_{\alpha_2}$ must belong to $A_{\alpha_1}$ and, thus, $c(\theta) \geq \alpha_1$, so that there exists no $\theta \in A_{\alpha_2}$ such that $c(\theta) = \alpha_2 < \alpha_1$, contradictory to the definition of $A_{\alpha_2}$.   

\bigskip
\noindent
(ii) Each $\widehat{\Theta}$ can be classified into 4 distinct types, denoted by (a) $\widehat{\Theta}_{\bar{LU}}$ where $\theta_L \not \in \widehat{\Theta}$ and $\theta_U \not \in \widehat{\Theta}$, (b) $\widehat{\Theta}_{LU}$ where $\theta_L \in \widehat{\Theta}$ and $\theta_U \in \widehat{\Theta}$ and, thus, $\theta \in \widehat{\Theta}_{LU}$, (c) $\widehat{\Theta}_{L}$ where $\theta_L \in \widehat{\Theta}$ and $\theta_U \not \in \widehat{\Theta}$, (d) $\widehat{\Theta}_{U}$ where $\theta_L \not \in \widehat{\Theta}$ and $\theta_U \in \widehat{\Theta}$. Type (c) can be further classified into (c.1) $\widehat{\Theta}_{L1}$ where $\theta \in \widehat{\Theta}_{L1}$ and (c.2) $\widehat{\Theta}_{L2}$ where $\theta \not \in \widehat{\Theta}_{L2}$, i.e. depending on whether or not $\theta$ appears in $\widehat{\Theta}$.  Similarly, type (d) into (d.1) $\widehat{\Theta}_{U1}$ where $\theta \in \widehat{\Theta}_{U1}$ and (d.2) $\widehat{\Theta}_{U2}$ where $\theta \not \in \widehat{\Theta}_{U2}$. We have
\begin{align*}
c(\theta_L) = \mathrm{Pr}(\widehat{\Theta}_{LU}) +\mathrm{Pr}(\widehat{\Theta}_L) & = \mathrm{Pr}(\widehat{\Theta}_{LU}) +\mathrm{Pr}(\widehat{\Theta}_{L1}) + \mathrm{Pr}(\widehat{\Theta}_{L2}) \\
c(\theta_U) = \mathrm{Pr}(\widehat{\Theta}_{LU}) +\mathrm{Pr}(\widehat{\Theta}_U) & = \mathrm{Pr}(\widehat{\Theta}_{LU}) +\mathrm{Pr}(\widehat{\Theta}_{U1}) + \mathrm{Pr}(r\widehat{\Theta}_{U2}) \\
c(\theta) & \geq \mathrm{Pr}(\widehat{\Theta}_{LU}) +\mathrm{Pr}(\widehat{\Theta}_{L1}) +\mathrm{Pr}(\widehat{\Theta}_{U1}).
\end{align*}
Thus, if $\mathrm{Pr}(\widehat{\Theta}_{U1}) \geq \mathrm{Pr}(\widehat{\Theta}_{L2})$, then $c(\theta) \geq c(\theta_L)$, or if $\mathrm{Pr}(\widehat{\Theta}_{U1}) \leq \mathrm{Pr}(\widehat{\Theta}_{L2})$, then $\mathrm{Pr}(\widehat{\Theta}_{L1}) \geq \mathrm{Pr}(\widehat{\Theta}_{U2})$ since $c(\theta_L) = c(\theta_U)$, such that $c(\theta) \geq c(\theta_U)$. Similarly on comparison between $\mathrm{Pr}(\widehat{\Theta}_{L1})$ and $\mathrm{Pr}(\widehat{\Theta}_{U2})$.  $\square$

\subsection{Proof of Theorem 2}

\noindent
Take any initial level-$\alpha_1$ corroboration set $A_{\alpha_1} = [L_{\alpha_1}, U_{\alpha_1}]$. Without losing generality, one of the end points must have corroboration $\alpha_1$ by Theorem 1.i; suppose $c(L_{\alpha_1}) \geq c(U_{\alpha_1}) = \alpha_1$. By definition $c(\theta) \geq \alpha_1$ for all $\theta\in A_{\alpha_1}$. If $c(\theta) = c(L_{\alpha_1})$ for all $L_{\alpha_1}< \theta < U_{\alpha_1}$, then $\theta^{\max} = L_{\alpha_1}$, since $c(\theta) < \alpha_1 \leq c(L_{\alpha_1})$ for any $\theta \not \in A_{\alpha_1}$. Otherwise, there exists $L_{\alpha_1}< \theta < U_{\alpha_1}$, where $c(\theta) = \alpha_2 > c(L_{\alpha_1}) \geq \alpha_1$, and the corresponding level-$\alpha_2$ corroboration set, denoted by $A_{\alpha_2} = [L_{\alpha_2}, U_{\alpha_2}]$. By Theorem 1.i, we have $[L_{\alpha_2}, U_{\alpha_2}] \subset [L_{\alpha_1}, U_{\alpha_1}]$. Since $\alpha\leq 1$, iteration of the argument must terminate at some maximum level-$\alpha$. $\square$

\subsection{Proof of Lemma 1}

\noindent
Let $\delta(\theta; \widehat{\psi}_n) =1$ if $\theta \in \mbox{Int}(\widehat{\Theta}) = (\widehat{L}_n, \widehat{U}_n)$, and $0$ otherwise, where $\widehat{\psi}_n$ is the MLE. Without losing generality, for any $\theta = U_0 - \epsilon$, where $0< 2 \epsilon < U_0 -L_0$, we have $\delta(\theta; \widehat{\psi}_n) =1$ if $|\widehat{U}_n - U_0| < \epsilon$ and $|\widehat{L}_n - L_0| < \epsilon$, the probability of which tends to 1, since $\widehat{\psi}_n \stackrel{\mbox{Pr}}{\rightarrow} \psi_0$. Thus, $\delta(\theta; \widehat{\psi}_n) \stackrel{\mbox{Pr}}{\rightarrow} 1$, i.e. $\bar{c}_0(\theta) = 1$ and $\theta \in \bar{A}^{\max}$. Similarly, it can be shown that $\bar{c}_0(\theta) = 1$, for $\theta\not \in \Theta_0$, i.e. $\theta \not \in \bar{A}^{\max}$. $\square$

\subsection{Proof of Theorem 3}

\noindent
By the general form of Slutsky's Theorem (e.g. Theorem 7.1, Kapadia \emph{et al.}, 2005), we have $\bar{c}(\theta; \widehat{\psi}_n) \stackrel{\mbox{Pr}}{\rightarrow} \bar{c}(\theta; \psi_0)$, since $\widehat{\psi}_n \stackrel{\mbox{Pr}}{\rightarrow} {\psi}_0$ and $\bar{c}(\theta; \psi)$ is a bounded for all $\psi$. Thus, if $\theta \in (L_0, U_0)$, such that $\bar{c}(\theta; \psi_0) = 1$ by Lemma 1, we have $\bar{c}(\theta; \widehat{\psi}_n) \stackrel{\mbox{Pr}}{\rightarrow} \bar{c}(\theta; \psi_0) =1$, meaning $\lim_n \mbox{Pr}(\theta \in \widehat{A}_n^{\max}) =1$. Similarly, it can be shown that  $\lim_n \mbox{Pr}(\theta \in \widehat{A}_n^{\max}) =0$, for $\theta\not \in \Theta_0$. $\square$

\end{document}